\documentstyle[12pt]{article}

\begin{document}
\title{Quantum Hall Effect Wave Functions
as Cyclic Representations of $U_q(sl(2))$}
\author{\"{O}mer F. DAYI\thanks{E-mail: dayi@mam.gov.tr.}\\
\small \it T\"{U}B\.{I}TAK--Research Institute for Basic Sciences,\\
\small \it P.O.Box 6, 81220 \c{C}engelk\"{o}y--Istanbul, Turkey. } 
\date{}
\maketitle
\begin{abstract}
Quantum Hall effect wave functions
corresponding to the filling factors
$\mbox{1/2p+1,}\ \mbox{2/2p+1,}\ \cdots ,\ \mbox{2p/2p+1,}\ 1,$
are shown to form a basis
of   irreducible cyclic representation
of the quantum algebra $U_q(sl(2))$ at $q^{2p+1}=1.$ 
Thus, the wave functions $\Psi_{P/Q}$ possessing  filling
factors $P/Q<1$ where $Q$ is odd and $P,\ Q$ are
relatively prime integers are classified in terms
of $U_q(sl(2)).$ 

\vspace{2cm}

\end{abstract}

\begin{flushleft}
RIBS--PH--TH--6/97

q-alg/9704010 
\end{flushleft}

\vspace{2cm}

\pagebreak

\vspace{1cm}
\noindent
{\large \bf 1. Introduction:}
\vspace{.5cm}

\noindent
Microscopic theory of the fractional
quantum Hall effect (QHE) is not well established.
Its theoretical understanding mostly is due to
trial wave functions \cite{qhe}.
For filling factors $1/m$ where $m$ is an odd integer, 
trial wave functions were given by 
Laughlin\cite{la} .
Trial wave functions for the other filling factors
$\nu =P/Q<1,$ where $P,\ Q$ are relatively prime integers
and $Q$ is odd,
were constructed in terms of some hierarchy
schemes\cite{owf}--\cite{jwz} 
where they were obtained from a parent
state which is a full filled Landau level or
a Laughlin wave function.
However, general properties of
the QHE should be independent of the explicit form of
trial wave functions, but depend on their universal
features as their  orthogonality.

We utilize orthogonality of QHE states for
different filling factors, independent of their
explicit form, to
show that they
can be classified
as irreducible cyclic representations
of $U_q(sl(2))$ at  roots of unity.
In our scheme,  states 
corresponding to 
filling factors  possessing a common denominator 
are in the same representation.

Although $U_q(sl(2))$ structures were found in 
the Hofstadter problem\cite{hof},
in the Landau problem\cite{lan},
for Laughlin wave functions\cite{laf}
and in the QHE\cite{uqh},
the approach presented here does not have any relation to 
them: {\it i)} As far as we deal with flat surfaces,
in all of the previous works
generators of the deformed algebra
were constructed in terms of magnetic translations.
The construction presented here cannot be
written in terms of magnetic transformations.
{\it ii)} Here, wave functions possessing different
filling factors which have a common denominator are treated
on the same footing. However, in the other works
only one state is considered and the theories
were built on them without mixing different states
with different parameters which  correspond to
filling factors in the QHE case.              

First, we show explicitly that the wave functions
corresponding to the filling factors $\nu = 1/3,\ 2/3, 1,$ 
can be considered as basis of
cyclic irreducible representation
of the quantum algebra $U_q(sl_2)$ at $q^3=1.$  
Then, the general case is studied.
Conclusions are presented in the last section.

\vspace{1cm}
\noindent
{\large \bf 2. Cyclic Representation of $U_q(sl(2)):$}
\vspace{.5cm}
\noindent

The deformed algebra $U_q(sl(2))$  
\begin{eqnarray}
[E_+,E_-] & = & \frac{K-K^{-1}}{q-q^{-1}},  \nonumber \\
KE_{\pm}K^{-1} & = & q^{\pm 2}E_{\pm}. \label{alg}
\end{eqnarray}
at roots of unity i.e. $q^{2p+1}=1,$ $p$ a positive integer, 
has a finite dimensional irreducible representation 
which has no classical finite dimensional analog.
This is the cyclic representation whose dimension 
is $2p+1$\cite{rep}. Cyclic means that there 
are no heighest or lowest weight states in the spectrum.
i.e. $E_+|\cdots >\neq 0$ and $E_-|\cdots >\neq 0$
for any state.

When $q^{2p+1}=1$ irreducible cyclic representation
of $U_q(sl(2))$
can be written in some basis
$\{ v_0,v_1, \cdots, v_{2p} \}$ as
\begin{eqnarray}
Kv_m & = & \lambda q^{-2m} v_m, \nonumber \\
E_+v_m & = & g_m v_{m+1} , \label{rep} \\
E_-v_m & = & f_m v_{m-1} , \nonumber
\end{eqnarray}
where $m=0,\cdots ,2p,$ and we defined $v_0\equiv v_{2p+1},\
v_{-1}\equiv v_{2p}.$ $\lambda,$ $g_m$ , and
$f_m$ are some complex constants which are nonzero and
in the case of requesting 
that the representation in
unitary, we should restrict their values
such that
\begin{equation}
\label{uni}
K^\dagger =K^{-1};\   E^\dagger_-=E_+.
\end{equation}

Although, for the purposes of this work there is no
need of discussing in detail neither how  unitary
representations arise in the general framework
nor  values of  Casimir  operators,
let us denote that there are three independent
Casimir operators of $U_q(sl(2))$ at $q^{2p+1}=1$:
$K^{2p+1},$ $E_+^{2p+1}$ and $E_-^{2p+1}.$  

\vspace{1cm}
\noindent
{\large \bf 3. Classification of $\nu= 1,\ 1/3,\ 2/3$
States:}
\vspace{.5cm}

\noindent
When $N$ particles (electrons) move on a plane in a
perpendicular magnetic field we may consider
the wave functions\cite{la,gi}
\begin{eqnarray}
\psi_{1} (z_1,\cdots ,z_N ) & = & {\cal N}_1 
e^{-\frac{1}{2} \sum_{k=1}^N|z_k|^2}\prod_{i<j}^N 
(z_i-z_j), \label{w1} \\
\psi_{1/3}(z_1,\cdots ,z_N )
& = & {\cal N}_2  
e^{-\frac{1}{2} \sum_{k=1}^N|z_k|^2}\prod_{i<j}^N 
(z_i-z_j)^3, \label{w2}  \\
\psi_{2/3}(z_1,\cdots ,z_N ) 
& = & {\cal N}_3 \int d^2z_{N+1}\cdots d^2z_{N+M}
e^{-\frac{1}{2} \sum_{K=1}^{N+M}|z_K|^2} \nonumber \\
 & & \prod_{l<n}^M (\bar{z}_{N+l}-
\bar{z}_{N+n})^3  
\prod_{I<J}^{N+M} (z_I-z_J), \label{w3}
\end{eqnarray}                         
which possess the following values of the angular momentum $L,$
\begin{eqnarray}
L [ \psi_{1}(z_1,\cdots ,z_N )]
& = & \frac{N(N-1)}{2} , \nonumber \\
L [ \psi_{1/3}(z_1,\cdots ,z_N )]
& = & 3\frac{N(N-1)}{2} ,         \nonumber  \\
L [ \psi_{2/3}(z_1,\cdots ,z_N )]
& = & \frac{(N+M)(N+M-1)}{2} -3\frac{M(M-1)}{2}. \nonumber 
\end{eqnarray}
It is supposed that $N$ is large
and we take $M=N/2.$

Filling factors of the $N$ particle
states are given in the 
thermodynamical limit as
\begin{equation}
\nu \equiv \lim_{N\rightarrow \infty}\frac{N(N-1)}{2L}.
\end{equation}
Hence, filling factors of
the wave functions (\ref{w1})--(\ref{w3}) are
\begin{equation}
\nu (\psi_1) =1,\ \nu (\psi_{1/3}) =1/3,\  
\nu(\psi_{2/3}) =2/3.
\end{equation}

Indeed, (\ref{w1}) is the wave function when the lowest
Landau level is fully filled  and (\ref{w2})--(\ref{w3})
are the trial wave functions which describe the QHE
at the filling  factors $1/3,\ 2/3.$

By making use of
\begin{equation}
\int d^2z e^{-|z|^2} \bar{z}^mz^n =\delta_{m,n},
\end{equation}
one can observe that
the wave functions which possess different angular momentum
values are orthogonal. Moreover, by
choosing the normalization constants ${\cal N}_a$ appropriately 
the wave functions 
(\ref{w1})--(\ref{w3}) can be taken  to satisfy ($N >2$)
\begin{equation}
\label{spr}
(\psi_\sigma , \psi_{\rho})\equiv
\int d^2z_1\cdots d^2z_N \bar{\psi}_\sigma (z_1\cdots ,z_N) 
\psi_{\rho} (z_1\cdots ,z_N) =\delta_{\sigma ,\rho},
\end{equation}
where $\sigma ,\rho =1,1/3,2/3.$

If  $\hat{\nu}$ denotes the first quantized operator
corresponding to the filling factor $\nu ,$
one can construct the physical operator
\begin{equation}
\label{dk}
\hat{k}\equiv e^{2\pi i \hat{\nu}},
\end{equation}
which will be shown to play the main role  in
classifying QHE wave functions in terms of $U_q(sl(2))$
at roots of unity.

In a second quantized theory, operators
corresponding to  physical operators of
the first quantization will be given in terms of 
states spanning the related field theory.
Let us deal with the states,
corresponding to (\ref{w1})--(\ref{w3}),
\begin{equation}
\label{3s}
|\sigma > \equiv  \int d^2z_{1}\cdots d^2z_N
e^{- \frac{1}{2}
\sum_{k=1}^N|z_k|^2}
\psi_\sigma(z_1,\cdots ,z_N)
|z_1,\cdots ,z_N>,
\end{equation}
where
\begin{equation}
\label{ve}
|z_1,\cdots ,z_N> = \frac{1}{\sqrt{N!}}
\varphi^\dagger (z_1) \cdots 
\varphi^\dagger (z_N) |0>.
\end{equation}
The fermionic operators $\varphi (z),\ \varphi^\dagger (z)$ 
satisfy the anticommutation relation
\[
\{\varphi^\dagger (z), \varphi (z^\prime )\} =e^{z^\prime \bar{z}}.
\]
The states (\ref{3s}) are orthonormal:
\begin{equation}
\label{sro}
<\sigma |\rho >=\delta_{\sigma ,\rho}.
\end{equation}

The second quantized operator 
\begin{equation}
\label{k12}
k=e^{2\pi i}|1><1 | +
e^{2\pi i/3} |1/3><1/3 |
+e^{4\pi i/3}|2/3><2/3 |
\end{equation}
corresponds to the first quantized physical
operator (\ref{dk}).
In terms of the vector
$(|1>,\ |1/3>,\ |2/3>)$ 
and the scalar product defined in (\ref{sro}),
one can obtain the representation
\begin{equation}
\label{k3}
k=\left(
\begin{array}{ccc}
1 & 0 & 0 \\
0 & \tilde{q} & 0 \\
0 & 0 & \tilde{q}^2
\end{array}
\right) ,
\end{equation}
where $\tilde{q}=\exp (2\pi i/3)$ i.e.
\[
\tilde{q}^3=1.
\]
Moreover, we can construct the operators
\begin{eqnarray}
e_+ & = &
a_1 |1/3><1|+a_2 |2/3><1/3 |+a_3 |1><2/3 |, \label{ep12} \\
e_- & = &
b_1 |1><1/3 |+b_2 |1/3><2/3 |+b_3|2/3><1 |,  \label{em12}
\end{eqnarray}
whose representations are
\begin{equation}
\label{e3}
e_+ =\left(
\begin{array}{ccc}
0 & a_1 & 0 \\
0 & 0 & a_2 \\
a_3 & 0 & 0
\end{array}
\right),\
e_- =\left(
\begin{array}{ccc}
0 & 0 & b_3 \\
b_1 & 0 & 0 \\
0 & b_2 & 0
\end{array}
\right)  .
\end{equation}
One can show that (\ref{k3}) and (\ref{e3}) realize
the $U_q(sl(2))$ algebra
\begin{equation}
\label{ans}
[e_+,e_-] =\frac{k-{k}^{-1}}{\tilde{q}
-\tilde{q}^{-1}} ,\
ke_{\pm}{k}^{-1}=\tilde{q}^{\pm 2}e_{\pm},
\end{equation}
if the coefficients satisfy 
\begin{eqnarray}
a_1b_1 -a_3b_3 & = & 0 , \nonumber \\
a_2b_2 -a_1b_1 & = & 1 , \label{ceq} \\
a_3b_3 -a_2b_2 & = & -1.  \nonumber
\end{eqnarray}
If one demands that the representation
(\ref{k3}), (\ref{e3})  is unitary:
\begin{equation}
{k}^{-1}={k}^\dagger ,\
{e_+}^\dagger =e_-,
\end{equation}
the coefficients $a_l,\ b_l$ should be
taken as
\begin{equation}
b_1=\bar{a}_1 ,\  b_2= \bar{a}_2 ,\ b_3=\bar{a}_3 .
\end{equation}
Then the conditions (\ref{ceq}) lead to
\begin{equation}
\label{cu}
|a_1|^2=|a_3|^2=|a_2|^2-1,\ a_l\neq 0 .
\end{equation}

Observe that the Casimir operators
\[
k^3= {\bf 1},\ e_+^3=a_1a_2a_3 {\bf 1},\  
e_-^3=b_1b_2b_3 {\bf 1}
\]
are proportional to identity.

An explicit realization  is presented in terms of the trial
wave functions (\ref{w1})--(\ref{w3}). However, the construction
depends only on the orthogonality of the states of
the QHE for different values of the filling factors and the
existence of the physical operator (\ref{dk}). This will be clarified
in the next section.

\vspace{1cm}
\noindent
{\large \bf 4. The General Case:}
\vspace{.5cm}

\noindent
QHE trial wave functions
in the standard hierarchy scheme
are given by\cite{owf},\cite{re}
\begin{equation}
\label{sth}
\psi_\nu (z_1,\cdots ,z_{N_0}) =\int
\prod_{\alpha =1}^{r}
\prod^{N_\alpha }_{i_\alpha =1}[d^2z^{(\alpha )}_{i_\alpha}]
e^{-\frac{1}{2}\sum_1^{N_0} |z_k|^2}
\prod_{\beta =0}^{r}
\prod_{i_\beta <j_\beta}^{N_\beta} 
(z^{(\beta )}_{i_\beta}-z^{(\beta )}_{j_\beta})^{a_\beta}
\prod_{i_{\beta+1},j_\beta =1}^{N_{\beta+1},N_\beta}
(z^{(\beta +1 )}_{i_{\beta +1}}-
z^{(\beta )}_{j_\beta})^{b_{\beta ,\beta +1}},
\end{equation}
where 
$z_{i_0}^{(0)} \equiv z_i .$
The measure 
$\prod [d^2z^{(\alpha )}_{i_\alpha}] $
depends on $a_\beta$ and 
$|z^{(\beta )}_{i_\beta}-z^{(\beta )}_{j_\beta}|,$ however the 
detailed form of it 
does not affect the filling factor $\nu =P/Q .$
$a_0$ is an odd positive integer,
$a_\alpha $ for $\alpha \neq 0$ are even integers
which can be positive or negative and
$b_{\beta +1, \beta}=\pm 1,$ except 
$b_{r,r+1}=0.$
By placing the $N_0$ electrons on a spherical
surface in a monopole magnetic field, one can
find that filling factor of (\ref{sth}) is given by
\begin{equation}
\label{gff}
\nu =\frac{1}{a_0 -\frac{1}{a_1- \frac{1}{\cdots -\frac{1}{a_r}}}}.
\end{equation}
Factors with negative powers may be replaced 
by complex--conjugate factors with positive powers
multiplied by some exponential factors.
Hence, (\ref{sth}) can equivalently be given as\cite{bw}
\begin{eqnarray}
\psi_\nu (z_1,\cdots ,z_{N_0})& = & \int
\prod_{\alpha =1}^{r}
\left[ \prod_{i_\alpha =1}^{N_\alpha} d^2z^{(\alpha )}_{i_\alpha}
\prod_{i_\alpha <j_\alpha}^{N_\alpha} 
|z^{(\alpha )}_{i_\alpha}
-z^{(\alpha )}_{j_\alpha}|^{2(-1)^{\alpha}\theta_\alpha} 
e^{-|q_\alpha |\sum_{i_\alpha}
|z^{(\alpha )}_{i_\alpha}|^2} \right] \nonumber \\
& & e^{-\frac{1}{2}\sum_1^{N_0} |z_{k}|^2}
\prod_{\beta =0}^{r} \prod_{i_\beta < j_\beta}^{N_\beta}
(\tilde{z}^{(\beta )}_{i_\beta} 
-\tilde{z}^{(\beta )}_{j_\beta})^{p_\beta}
\prod_{i_{\beta+1},j_\beta=1}^{N_{\beta+1},N_\beta}
(\bar{\tilde{z}}_{i_{\beta +1}}^{(\beta +1) }
-\tilde{z}^{(\beta )}_{j_\beta}), \label{bws}
\end{eqnarray}
where $\tilde{z}^{(\beta)}_{i_\beta} =z^{(\beta )}_{i_\beta}$ 
for $\beta =$ even
and $\tilde{z}^{(\beta)}_{i_\beta} =\bar{z}^{(\beta)}_{i_\beta}$ 
for $\beta =$ odd and
\begin{eqnarray}
\theta_0=0,  & \theta_r  = 
\frac{(-1)^r}{p_{r-1}-(-1)^r\theta_{r-1}}, \nonumber \\
q_0  =-1, & q_r  =   (-1)^{r+1}q_{r-1}\theta_r. \nonumber
\end{eqnarray}
Now, filling factor is
\begin{equation}
\label{ffbw}
\nu =\frac{1}{p_0 +\frac{1}{p_1+ \frac{1}{\cdots +\frac{1}{p_r}}}},
\end{equation}
where $p_0$ is odd and the other $p_i$ are even integers.

By generalizing the calculations of 
Laughlin given in Ref. \cite{qhe}
and making use of the scalar product defined in (\ref{spr}),
one can show that $\psi_\nu$ states are orthogonal\cite{re}.

To emphasize the second quantized character of our construction
let us introduce the states
\begin{equation}
\label{gb}
|i,\ p>_T  =  \int d^2z_{1}\cdots d^2z_{N_0}
e^{- \frac{1}{2}
\sum_{k=1}^{N_0}|z_k|^2}
\psi_{\frac{i}{2p+1}}(z_1,\cdots ,z_{N_0})
|z_1,\cdots ,z_{N_0}>,
\end{equation}
where $i=1,\cdots ,2p+1;\  p=1,2, \cdots ,$  
so that any filling factor $\nu =P/Q$ is considered.
We used the vectors (\ref{ve}) with $N$ replaced by $N_0.$
The subscript $T$ denotes the fact that trial wave
functions are used to give an explicit realization.

The states (\ref{gb}) are orthonormal:
\begin{equation}
\label{ort}
_T<i,\ p|j,\ p^\prime >_T=\delta_{i,j}\delta_{p,p^\prime}.
\end{equation}

We have shown that the states
$|i,\ p>_T $ are orthonormal by using the explicit form
of trial wave functions. However,
this should be a universal feature of QHE wave functions.
Then,
even if we do not know the explicit form, we can say that
exact states of the QHE  which we  indicate with $|i,\ p>,$
should be orthonormal:
\begin{equation}
\label{ort1}
<i,\ p|j,\ p^\prime >=\delta_{i,j}\delta_{p,p^\prime}.
\end{equation}
Indeed, in the following we will use this universal
property of QHE states without referring to any
trial wave function.

To generalize the construction given in Section 3,
let us deal with the states
\begin{equation}
\label{hik}
|1,\ p>,\ |2,\ p>,\cdots,\ |2p,\ p>,\ |2p+1,\ p>,
\end{equation}
corresponding to the filling factors, respectively,
\begin{equation}
\label{gnu}
\nu =\frac{1}{2p+1}, \frac{2}{2p+1},\cdots ,\frac{2p}{2p+1},1.
\end{equation}

Define the following
second quantized operators acting in the space
spanned by the states  (\ref{hik}),
\begin{eqnarray}
\tilde{K}
& = & \sum_{i=1}^{2p+1} q^i|i,\ p>
<i,\ p| ,\label {gc0} \\
\tilde{E}_+
& = & \sum_{i=1}^{2p+1} a_i|i,\ p>
<i+2,\ p| , \label{gc1} \\
\tilde{E}_-
&= & \sum_{i=1}^{2p+1} \bar{a}_i
|i+2,\ p><i,\ p| , \label{gc2}
\end{eqnarray}
where 
\begin{equation}
\label{root}
q^{2p+1} =1.
\end{equation}
To obtain the compact forms we adopted the definitions
\[
|2p+2,\ p>\equiv |1,\ p> ,\  |2p+3,\ p>\equiv |2,\ p>.
\]

By using the orthonormality condition (\ref{ort})
one observes that
inverse of $\tilde{K}$ is 
\begin{equation}
\label{k-1}
{\tilde{K}}^{-1} =  \sum_{i=1}^{2p+1} q^{-i}|i,\ p><i,\ p| =
{\tilde{K}}^\dagger .
\end{equation}

Let the coefficients $a_i$ 
are nonzero and satisfy
\begin{eqnarray*}
|a_{2p+1}|^2  - |a_{2p-1}|^2 & = & 0, \\
|a_{2p}|^2 -  |a_{2p-2}|^2 & = & -1, \\
|a_{l+2}|^2  -  |a_{l}|^2 & = &
\frac{q^{l+2}-q^{-l-2}}{q-q^{-1}}, \\
\end{eqnarray*}
where $l=-1,0,\cdots (2p-3);\ a_{-1}\equiv a_{2p},\
a_0\equiv a_{2p+1}.$ 
Then,
in terms of the basis $(|1,\ p>,\ \cdots ,\ |2p+1,\ p>)$
the operators
(\ref{gc0})--(\ref{gc2}) lead to a $(2p+1)$
dimensional unitary irreducible cyclic  representation
of $U_q(sl(2))$ at q satisfying (\ref{root}).

Note that the Casimir operators are proportional to unity as before:
$\tilde{K}^{2p+1}={\bf 1}$ and
$\tilde{E}_+^{2p+1}={\tilde{E}_-}^{\dagger\ 2p+1}=
\left( \prod_{i=1}^{2p+1}a_i\right) {\bf 1}.$

\vspace{1cm}
\noindent
{\large  \bf 5. Discussions:}
\vspace{.5cm}

\noindent
It is shown that QHE wave functions can be classified as
irreducible cyclic representations of
$U_q(sl(2))$ at roots of unity in a very natural way.
This naturalness follows from the fact that
the most significant physical quantity
of the QHE $\nu =P/Q $ 
fits very well with the integer ($m$ in (\ref{rep}))
characterizing irreducible cyclic representations of 
$U_q(sl(2)).$  
Obviously,  any set of orthogonal
states possessing a quantum number
which permits a partition of 
unity like $\nu ,$ 
\[
\sum_{i=1}^{2p+1}\frac{ \nu (|i,\ 2p+1>)}{p+1} =1 ,
\]
can be classified 
as irreducible cyclic  representation
of $U_q(sl(2))$ at a root of unity.

How one can utilize the proposed classification of
the QHE to calculate some physical quantities?
Here, one of the most significant physical quantities is
the partition function which may be obtained if the
Green function in the space defined by
$U_q(sl(2))$ at roots of unity with cyclic representation
is available. In Ref. \cite{ad} Green function in the space
defined by the q--deformed group $SU_q(2)/U(1)$
for $q$ not a root of unity is obtained
without referring to explicit forms of the representations
but depending only on their general features.
We hope that a similar calculation can be used
in our case. Then, we can obtain Green function and
in terms of that the related
partition function which
may give some hints about its
physical interpretation which is not clear at the moment.

\vspace{.5cm}

\noindent
{\bf Acknowlegment:}

\vspace{.3cm}

\noindent
I would like to thank I.H. Duru for fruitful discussions.

\pagebreak


\begin{thebibliography}{99}
\bibitem{qhe}{\it The Quantum Hall Effect,} edited by 
R. E. Prange and S. M. Girvin (Springer-Verlag, New York, 1987).
\bibitem{la}R. B. Laughlin, Phys. Rev. Lett. {\bf 50 } (1983) 1395.
\bibitem{owf}F. D. M. Haldane, Phys. Rev. Lett. {\bf 50 }  (1983) 605;\\
B. I. Halperin, Phys. Rev. Lett. {\bf 52 } (1984) 1583.
\bibitem{jwz}J. K. Jain, Phys. Rev. Lett. {\bf 63}  (1989) 199;\\
X. G. Wen and A. Zee, Phys. Rev. B {\bf 46 } (1992) 2290.
\bibitem{hof}P.B. Wiegmann and A.V. Zabrodin, Nucl. Phys. B
{\bf 422} (1994) 495;\\
L.D. Faddeev and R.M. Kashaev, Comm. Math. Phys. {\bf 169}
(1995) 181.
\bibitem{lan} H-T. Sato, Mod. Phys. Lett. A {\bf 9 } (1994) 451;\\
C-L. Ho, Mod. Phys. Lett. A {\bf 10 } (1995) 451
\bibitem{laf}N. Aizawa, S. Sachse and H-T. Sato,
Mod. Phys. Lett. A {\bf10 } (1995) 853;\\
M. Alimohammadi and H. Mohseni Sadjadi, J. Phys. A {\bf 29} (1996) 5551.
\bibitem{uqh}I. I. Kogan, Int. J. Mod. Phys. A {\bf 9 } (1994) 3887; \\
H-T. Sato, Mod. Phys. Lett. A {\bf 20 } (1994) 1819;\\
M. Alimohammadi and A. Shafei Deh Abad, J. Phys. A {\bf 29 }
(1996) 559; \\
G-H. Chen, L-M. Kuang and M-L. Ge, Phys. Rev. B {\bf 53 } (1996) 9540.
\bibitem{rep}G. Lusztig, Adv. Math. {\bf 70 } (1988) 237;\\
M. Rosso, Commun. Math. Phys. {\bf 117 } (1988) 581;\\
P. Roche and D. Arnaudon, Lett. Math. Phys. {\bf 17 } (1989) 295;\\
C. de Concini and  V. G. Kac, {\it Representations of Quantum 
Groups at Root of Unity,} Progress in Mathematics Vol. 92 
(Birkh\"{a}use, Boston, 1990); \\
L. C. Biedenharn and M. A. Lohe, {\it Quantum Group Symmetry and
q--Tensor Algebras,} (World Scientific, Singapore, 1995).
\bibitem{gi}S. M. Girvin, Phys. Rev. B {\bf 29 }  (1984) 6012.
\bibitem{re}N. Read, Phys. Rev. Lett. {\bf 65} (1990) 1502.
\bibitem{bw}B. Blok and X. G. Wen, 
Phys. Rev. B {\bf 43 } (1991) 8337.
\bibitem{ad}H. Ahmedov and I.H. Duru, ``Green Function
on the q--Symmetric Space $SU_q(2)/U(1)$", RIBS-PH-4/97,
q-alg/9703032.
\end{thebibliography}
\end{document}